\documentclass[english,russian]{amsart}
\usepackage{amssymb,amsmath}
\newtheorem{thm}{Theorem}[section]

\newtheorem{lm}[thm]{Lemma}

\newcommand{\sudda}[1]{}

\begin{document}

\title{Assosymmetric algebras under Jordan product  }

\author{Askar Dzhumadil'daev}
\address
{Kazakh-British University, Tole bi 59, Almaty,050000,Kazakhstan}
\email{dzhuma@hotmail.com} 



\maketitle

\begin{abstract}
We  prove that  assosymmetric algebras under Jordan product are  Lie triple. A Lie triple algebra is called special if it is isomorphic to a subalgebra of some plus-assosymmetric algebra. We establish that Glennie identitiy  is valid for special Lie triple algebras, but not 
 for all Lie triple algebras.
\end{abstract}

An associative algebra under Jordan commutator satisfies three remarkable polynomial identities: commutativity, Jordan identity and Glennie identity. 
These identities are independent. It is not clear whether these identities form base of identities for plus-associative algebras \cite{Shestakov}.

In our paper we establish analogue  of these results for assosymmetric algebras. Assosymetric algebras play in  theory of Lie triple algebras the same role as a role of associative algebras in  theory of Jordan algebras. We prove that any assosymmetric algebra under Jordan product satisfies three polynomial identities: commutativity, Lie triple identity and Glennie identity. In a class of plus-assosymmetric algebras these identities are independent. 
 Lie triple identity is a consequence of Jordan identity, but these identities are not equivalent. 

To formulate our results we need to introduce some notations and remind some definitions. Let $K$ be a field of characteristic $p\ge 0.$  We will suppose that $p\ne 2,3, $ if otherwise is not stated.
 For algebra $A$ notation $A=(A,\cdot)$ will mean that a vector space $A$ over field $K$ has multiplication $\cdot,$ i.e., 
$A$ is endowed by  bilinear map $(a,b)\mapsto a\cdot b.$ For an algebra $A$ with multiplication $\cdot$ its {\it plus-algebra} is defined as algebra 
$A^{(+)}=(A,\star),$ where $a\star b=a\cdot b+b\cdot a$ is Jordan commutator. Similarly a {\it minus-algebra} $A^{(-)}=(A, [\;,\;])$  of $A$ is defined as  algebra with vector space $A$ and multiplication given by Lie commutator $[a,b]=a\cdot b-b\cdot a.$
For a class of algebras $\mathcal C$ denote by ${\mathcal C}^{(\pm )}$ class of algebras of a form $A^{(\pm)},$ where $A\in \mathcal C.$
 
Let   $F=K\langle t_1,t_2,\ldots\rangle$ be an algebra of non-associative non-commutative polynomials
with countable number of  generators  $t_1,t_2,\ldots$ . Sometimes we will restrict the number of generators and we consider  
$F$ as an absolute  free algebra with generators $t_1,t_2,\ldots,t_k.$ The algebra $F$ sometimes is called absolute free algebra, sometimes free magmatic algebra. 
For an algebra $A=(A,\cdot)$ with multiplication $\cdot$ and a polynomial $f=f(t_1,\ldots,t_k)\in F$ say that $f=0$ is a {\it polynomial identity} of $A$ if 
$f(a_1,\ldots,a_k)=0$ for any substitution $t_i:=a_i\in A,$ $1\le u\le k.$ Here we calculate $f(a_1,\ldots,a_k)$ in terms of multiplication $\cdot.$ 
About polynomial identities and varieties  of non-associative algebras see for example \cite{Shestakov}.

Let $[t_1,t_2]=t_1t_2-t_2t_1$ be Lie commutator, 
$t_1\star t_2=t_1t_2+t_2t_1$ Jordan commutator and 
$$(t_1,t_2,t_3)=t_1(t_2t_3)-(t_1t_2)t_3$$
associator.  Associator of plus-algebra is denoted by $\langle a,b,c\rangle ,$
$$\langle t_1,t_2,t_3\rangle =t_1\star(t_2\star t_3)-(t_1\star t_2)\star t_3.$$
Let 
$$lsym(t_1,t_2,t_3)=(t_1,t_2,t_3)-(t_1,t_3,t_2),$$
$$rsym(t_1,t_2,t_3)=(t_1,t_2,t_3)-(t_1,t_3,t_2),$$
be left-symmetric and right-symmetric identities.  An algebra with identities $lsym=0$ and $rsym=0$  is called {\it assosymmetric.} 
So, an algebra $A=(A,\cdot)$  is assosymmetric if 
$$a\cdot [b,c]=(a\cdot b)\cdot c-(a\cdot c)\cdot b,$$
$$[a,b]\cdot c=a\cdot(b\cdot c)-b\cdot(a\cdot c).$$
for any $a,b,c\in A.$ Recall that {\it nucleus} of algebra $A$  is a subspace of elements $x\in A$ such that 
$(x,a,b)=(a,x,b)=(a,b,x)=0$ for any $a,b\in A.$ A nucleus of assosymmetric algebras contains  an ideal generated by commutator $[A,A]$ (\cite{Boers}).
In our paper we will use the following properties of assosymmetric algebra $A,$ 
$$([a,b],c,d)=(a,[b,c],d)=(a,b,[c,d])=0,$$
$$[a,b]\cdot(c,d,e)=(a,b,c)\cdot [d,e]=0,$$
for any $a,b,c,d,e\in A.$ Assosymmetric algebras was studied in 
 (\cite{Boers},  \cite{Hentzel},  \cite{Hentzel1},   \cite{Kleinfeld},    \cite{Osborn2}, \cite{Petersson}, \cite{Sitaram}). Free base of assosymmetric algebras was found in  \cite{Hentzel}.
 
 Denote by $F_{\mathcal C}^{\pm }(t_1,\ldots,t_k)$ subalgebra of $F_{\mathcal C}(t_1,\ldots,t_k)^{(\pm )}$ 
generated by elements $t_1,\ldots,t_k.$ Then $F_{{\mathcal A}ss}^-(t_1,\ldots,t_k)$ is isomorphic to free Lie algebra generated by $t_1,\ldots,t_k.$ In that time 
$F_{\mathcal{A}ss}^+$ is not free.   For element $X\in F_{\mathcal C}(t_1,\ldots,t_k)$ say that $X$ is {\it Lie element} in the class $\mathcal C$ if $X\in  F^-_{\mathcal C}(t_1,\ldots,t_k).$ Similarly, $X\in F_{\mathcal C}(t_1,\ldots,t_k)$ is {\it Jordan element} in the class 
 $\mathcal C$ if $X\in  F^+_{\mathcal C}(t_1,\ldots,t_k).$ If ${\mathcal C}=\mathcal{A}ss$ then Jordan element is called {\it $j$-element.} 
 If ${\mathcal C}={\mathcal A}ssym$ then its Jordan element is called {\it  $ja$-element. } If ${\mathcal C}={\mathcal A}lt$ is a class of alternating algebras, then 
its  Jordan element is called as {\it $jalt$-element.}

 Let 
 $$jor(t_1,t_2,t_3)=(t_1,t_2,t_1^2),$$
$$wjor(t_1,t_2,t_3,t_4)=(t_2,t_1,t_3t_4)+(t_3,t_1,t_4t_2)+(t_4,t_1,t_2t_3)$$
be Jordan and multilinear Jordan polynomials. 
If $p\ne 2,3,$ then identities $jor=0$ and $wjor=0$ are equivalent. If $p=3,$ then the identity $wjor=0$ is a consequence of identity $jor=0,$ but converse is not true. Recall that Jordan algebras are defined as commutative algebras with identity $jor=0.$  Associative algebras under Jordan product satisfies  the Jordan identity of degree $4,$
$$\langle a,b,a^2\rangle =0$$
and one more identity of degree $8,$ so called Glennie identity \cite{Glennie}. 
We give Shestakov's construction (\cite{McCrimmon}, \cite{Sverchkov1}) of Glennie identity. Let 
$$shest(t_1,t_2,t_3)=-3\langle t_1,t_3,t_2\rangle  \star (\langle t_1,t_1,t_2^2\rangle-\langle t_1,t_1,t_2\rangle\star t_2)-2\langle t_1,\langle t_1,\langle t_1,t_3,t_2\rangle,t_2\rangle,t_2\rangle$$ 
be Shestakov polynomial and
$$glen(t_1,t_2,t_3)=shest(t_1,t_2,t_3\star t_3)-2t_3\star shest(t_1,t_2,t_3)$$
be Glennie polynomial. 
Shestakov has established that the derivation $D_{([a,b]\star [a,b])\star [a,b]}\in Der\,F_{{\mathcal A}ss}$ is well defined on Jordan subspace:
for any $a,b,c\in F_{Ass},$ the element 
$$D=D(a,b,c)\stackrel{def}{=}[([a,b]\star [a,b])\star [a,b],c]$$ 
is $j$-element. He got exact construction of the element $D(a,b,c)$ as $j$-element,
\begin{equation} \label{shestakov}
D(a,b,c)=shest(a,b,c).
\end{equation}
Then Glennie identity is equivalent to the Leibniz condition for derivation 
 $$D_{([a,b]\star [a,b])\star [a,b]}\in Der\,F_{{\mathcal A}ss}^{(+)}.$$

In our paper we prove that relation (\ref{shestakov}) holds not only for associative algebras, but also for assosymmetric algebras.
In other words, the element $D(a,b,c)$ is not only $j$-element, but also $ja$-element. 
Note that relation (\ref{shestakov}) fails for a class alternative algebras. It does not mean that $D(a,b,c)$ is not $jalt$-element. We would be surprised very much if it is so. It just means that in alternative case Jordan polynomial corresponding to the element $D(a,b,c)$ differs from Shestakov polynomial. 

Let 
$$lietriple(t_1,t_2,t_3)=(t_1,t_2^2,t_3)-t_2\star (t_1,t_2,t_3)$$
be Lie triple identity. 
A commutative algebra with identity $lietriple=0$ is called {\it Lie triple} (\cite{JordanMatsushita} , \cite{Osborn2}, \cite{Petersson}.
So, $A$ is  Lie triple,  if for any $a,b\in A,$ 
$$(a,b^2,c)=2b(a,b,c).$$

Denote by $\mathcal{A}ss^{(+)}$ class of plus-associative algebras, i.e. class of algebras  $A^{(+)},$ where $A$ runs associative algebras. 
Similarly, $\mathcal{A}ssym^{(+)}$ is a class of assosymmetric algebras under Jordan product.

\begin{thm} \label{main1} Let $p\ne 2.$ Any plus-algebra of assosymmetric algebra $A$ is Lie triple, 
$$\langle a,b\star b,c\rangle =2\,b\star\langle a,b,c\rangle,\qquad \forall a,b,c\in A,$$
and satisfies  Glennie identity
$$glen(a,b,c)=shest(a,b,c\star c)-2c\star shest(a,b,c)=0,\qquad \forall a,b,c\in A.$$
The identities $[t_1,t_2]=0, lietrriple(t_1,t_2,t_3)=0$ and $glen(t_1,t_2,t_3)=0$ are independent. 
\end{thm}

\begin{thm} \label{q^2=1}  If $p\ne 2,3,$ then any identity of $\mathcal{A}ssym^{(+)}$  of degree $4$ follows from commutativity and 
Lie triple identities.

If $p=3,$ then Lie triple identity is not minimal for the class ${\mathcal A}ssym^{(+)}$. It satisfies the identity $wjor=0$ and 
any identity of degree $4$ for the class $\mathcal{A}ssym^{(+)}$ in case $p=3$ follows from commutativity one and the identity $wjor=0.$
\end{thm}

\section{\label{Associators} Associators, derivations and commutators  for assosymmetric algebras}

In this section we assume that $A=(A,\cdot)$ is assosymmetric if otherwise is not stated. Recall that 
 $a\star b=a\cdot b+b\cdot a,$ $[a,b]=a\cdot b-b\cdot a$  and $\langle a,b,c\rangle=a\star(b\star c)-(a\star b)\star c$ be associator of Jordan multiplication.
Here we give some preliminary results that we will use in proof of our theorems. 
In proof of our lemmas we  use results of \cite{Boers} . This paper needs  restriction $p\ne 2,3.$ 

\begin{lm} \label{one} For any algebra $A=(A,\cdot),$ 
$$\langle a,b,c\rangle -(a,b,c)+(c,b,a)=a\cdot(c\cdot b)-c\cdot(a\cdot b)-(b\cdot a)\cdot c+(b\cdot c)\cdot a.$$
\end{lm}

\begin{lm}\label{two} For any $a,b,c\in A,$
$$\langle a,b,c\rangle=[[a,c],b].$$
\end{lm}

{\bf Proof.} By Lemma \ref{one} and by assosymmetric rules
$$\langle a,b,c \rangle=a\cdot(c\cdot b)-c\cdot(a\cdot b)-(b\cdot a)\cdot c+(b\cdot c)\cdot a=$$
$$[a,c]\cdot b-b\cdot [a,c]=[[a,c],b].$$

\begin{lm} \label{Boers} 
If $A$ is assosymmetric, and $u=[x,y]$ is commutator element for some $x,y\in A,$ then 
$$(u,a,b)=(a,u,b)=(a,b,u)=0,$$ 
and 
$$u\cdot (a,b,c)=(a,b,c)\cdot u=0.$$ 
Moreover, if $a$ has a form $[x,y]\cdot z$ or $x\cdot [y,z],$ for some $x,y,z\in A,$ then also 
$$(u,a,b)=(a,u,b)=(a,b,u)=0.$$ 
\end{lm}

{\bf Proof.} Follows from results of \cite{Boers}.

\begin{lm}\label{der} Let  $ad\,a:A\rightarrow A$  be an adjoint map, $ad\,a(b)=[a,b]=a\cdot b-b\cdot a.$
Then 
$$ad\,a(b\cdot c)=ad\,a(b)\cdot c+b\cdot ad\,a(c)+(a,b,c).$$
for any $a,b,c\in A.$ If $a$ has a form $a=[x,y]$ or $a=[x,y]\cdot z,$ or $a=x\cdot [y,z],$ for some $x,y,z\in A,$ 
then $ad\,a\in Der\,A.$ 

In particular, $ad\,a$ is a derivation of minus-algebra $A^-=(A,[\;, \;]),$ 
$$ad\,a [b,c]=[ad\,a(b),c]+[a,ad\,a(c)].$$

The map $ad\,a$ is not derivation of plus-algebra $A^+=(A, \star ),$
$$ad\,a (b\star c)=ad\,a(b)\star c+b\star ad\,a(c)+2(a,b,c).$$
But $ad\,a$ is a derivation of plus-algebra, if $a\in [A,A].$ 
\end{lm}

{\bf Proof.} Follows from Lemma \ref{Boers}.

\begin{lm}\label{K3}
$$[(a,b,b), a]\star [[a,b],c]=0.$$
\end{lm}

{\bf Proof.} Let $ u=[[a,b],c].$ Then by Lemma \ref{Boers}
$$((a,b,b)\cdot a)\cdot u=(a,b,b)\cdot [a,u]+((a,b,b)\cdot u)\cdot a=0,$$
$$u\cdot(a\cdot (a,b,b))=[u,a]\cdot(a,b,b)+a\cdot(u\cdot(a,b,b))=0.$$

Therefore,
$$[(a,b,b), a]\star u=$$
$$[(a,b,b),a]\cdot u+u\cdot [(a,b,b),a]=$$
$$((a,b,b)\cdot a)\cdot u-(a\cdot (a,b,b))\cdot u+u\cdot ((a,b,b)\cdot a)-u\cdot(a\cdot(a,b,b))$$
$$-(a\cdot (a,b,b))\cdot u+u\cdot ((a,b,b)\cdot a)=$$

$$-[a\cdot (a,b,b),u]
-u\cdot (a\cdot (a,b,b)))
+((a,b,b)\cdot a)\cdot u
-[(a,b,b)\cdot a,u]=$$

$$-[a\cdot (a,b,b),u]-[(a,b,b)\cdot a,u]=$$
So, by Lemma \ref{der} and Lemma \ref{Boers} ,
$$[(a,b,b), a]\star u=$$
$$-[a,u]\cdot (a,b,b)-a\cdot [(a,b,b),u]
-[(a,b,b),u]\cdot a
-(a,b,b)\cdot [a,u]=$$
$$0.$$

\begin{lm}\label{five} Let $A=(A,\cdot)$ be assosymmetric algebra. For $a,b,c\in A$ set 
$$D=D(a,b,c)=[([a,b]\star [a,b])\star[a,b],c],$$
$$shest=shest(a,b,c)=-3\langle a,c,b\rangle  \star (\langle a,a,b^2\rangle-\langle a,a,b\rangle\star b)-2\langle a,\langle a,\langle a,c,b\rangle,b\rangle,b\rangle.$$ 
Then $D(a,b,c)=shest(a,b,c).$ In particular, $D(a,b,c)$ is Jordan element of assosymmetric algebra generated by $a,b,c.$
\end{lm}

{\bf Proof.}  Let us set 
$$x=[a,b], \quad u=[[a,b],c].$$
Then $u=[x,c].$ We see that $u$ and $x$ are commutator elements. 

By Lemma \ref{der}
$$D=[(x\star x)\star x,c]=$$
$$-2[c,(x\cdot x)\star x]=$$
$$-2[c,x\cdot x]\star x-2(x\cdot x)\star [c,x]-4\,(c,x\cdot x,x)=$$
$$-2([c,x]\cdot x)\star x-2(x\cdot[c,x])\star x-2(x\cdot x)\star [c,x]-2(c,x, x)\star x-4\,(c,x\cdot x,x)=$$
$$-2([c,x]\star x)\star x-2(x\cdot x)\star [c,x]-2(c,x, x)\star x-4\,(c,x\cdot x,x).$$
Since $x\in [A,A],$ by Lemma \ref{Boers}.
$$(c,x, x)=(c,x,x\cdot x)=0.$$
Therefore, 
$$D=2\,x\star(x\star u)+(x\star x)\star u.$$

Let us set
$$D_2= [[a,b^2],a]-[[a,b],a]\star b.$$
By Lemma \ref{der}
$$[a,b^2]=[a,b]\cdot b+b\cdot [a,b]+(a,b,b)=x\cdot b+b\cdot x+(a,b,b)$$
Therefore,
$$D_2=  [[a,b^2],a]-[[a,b],a]\star b=$$
$$[x\cdot b,a]+[b\cdot x,a]+[(a,b,b),a]-[x,a]\star b.$$
So, by Lemma \ref{der}
$$D_2=[x,a]\cdot b+x\cdot[b,a]-(a,x,b)
+[b,a]\cdot x+b\cdot[x,a]-(a,b,x)
+[(a,b,b),a]
-[x,a]\cdot b-b\cdot [x,a]=$$
$$-2x\cdot x+[(b,a,b),a]-2(a,b,x).$$

Then by Lemma \ref{two}
$$shest=-3[[a,b],c]\star ([[a,b^2],a]-[[a,b],a]\star b)-2 [[a,b],[[a,b],[[a,b],c]]]=$$ 
$$-3 D_2\star u-2[x,[x,u]].$$
Thus,
$$D-shest=D_1+3D_2\star u,$$
where we set 
$$D_1=2\,x\star(x\star u)+2[x,[x,u]]+(x\star x)\star u.$$
 
We have 
$$D_1=$$
$$2\,x\cdot(x\cdot u)+2\,x\cdot(u\cdot x)+2\,(x\cdot u)\cdot x+2\,(u\cdot x)\cdot x$$
$$+2\,x\cdot(x\cdot u)-2\,x\cdot(u\cdot x)-2\,(x\cdot u)\cdot x+2\,(u\cdot x)\cdot x$$
$$+(x\star x)\star u=$$
$$4\,x\cdot(x\cdot u)+4\,(u\cdot x)\cdot x+(x\star x)\star u.$$
Thus,
$$D_1+3D_2\star u=$$
$$4\,x\cdot(x\cdot u)+4\,(u\cdot x)\cdot x-6(x\cdot x)\cdot u-6\,u\cdot(x\cdot x)+3[(a,b,b),a]\star u-6(a,b,x)\star u+(x\star x)\star u=$$
$$4\,(x\cdot(x\cdot u)-(x\cdot x)\cdot u)+4\,((u\cdot x)\cdot x-u\cdot(x\cdot x))
-2(x\cdot x)\cdot u-2\,u\cdot(x\cdot x)+3[(a,b,b),a]\star u-6(a,b,x)\star u+(x\star x)\star u=$$
$$4\,(x,x,u)-4\,(u,x,x)-(x\star x)\star u+3[(a,b,b),a]\star u-6(a,b,x)\star u+(x\star x)\star u=$$
$$(3[(a,b,b),a]-6(a,b,x))\star u.$$
Since by Lemma \ref{Boers}, $(a,b,x)\star u=0,$ we see that 
$$D_1+3D_2\star u=3[(a,b,b),a]\star u.$$
It remains to use Lemma \ref{K3}, to obtain that 
$$D-shest=D_1+3\,D_2\star u=0.$$

\begin{lm}\label{4}
$\langle b,a,c\star d\rangle+\langle c,a,d\star b\rangle+\langle d,a, b\star c\rangle=-6[a,(b,c,d)]$
\end{lm}

{\bf Proof.} We have 
$$[a,b\cdot c]+[b,c\cdot a]+[c,a\cdot b]=(a,b,c)+(b,c,a)+(c,a,b)=3(a,b,c).$$
Hence,
$$[a,b\star c]+[b,c\star a]+[c,a\star b]=6(a,b,c).$$
Therefore, by Lemma \ref{two},
$$\langle b,a,c\star d\rangle+\langle c,a,d\star b\rangle+\langle d,a, b\star c\rangle=$$
$$[[b,c\star d],a]+[[c,d\star b],a]+[[d,b\star c],a]=$$
$$[[b,c\star d]+[c,d\star b]+[d,b\star c],a]=$$
$$-6[a,(b,c,d)]$$

\begin{lm}\label{5}
For any assosymmetric algebra $A$ a polylinear map 
$$A\times A\times A\times A\rightarrow A, \quad (a,b,c,d)\mapsto wjor(a,b,c,d)$$
is symmetric,
$$wjor(a_1,a_2,a_3,a_4)=wjor(a_{\sigma(1)},a_{\sigma(2)},a_{\sigma(3)},a_{\sigma(4)}),$$
for any permutation $\sigma\in S_4.$
\end{lm}

{\bf Proof.}  
If $\sigma\in S_4$ fixes the first element, $\sigma(1)=1,$ then 
$$wjor(a_1,a_2,a_3,a_4)=wjor(a_{\sigma(1)},a_{\sigma(2)},a_{\sigma(3)},a_{\sigma(4)}),$$
for any algebra $A$, not necessary assosymmetric.  The property
$$[a_1,(a_2,a_3,a_4)]=[a_2,(a_1,a_3,a_4)]$$
was  established in \cite{Boers}, p.14, the relation (vii) or (vii'). Thus, by Lemma~\ref{4} 
$wjor(a_1,a_2,a_3,a_4)$ is symmetric by all parameters.  

\section{Proof of Theorem \ref{main1}}

 By Lemma \ref{two}
$$\langle a,b\star b,c\rangle=2\,[[a,c],b^2]=$$ 
$$2\,[[a,c],b]\cdot b+2\,b\cdot[[a,c],b]=2 \langle a,b,c\rangle \cdot b+2\,b\cdot \langle a,b,c\rangle=2\langle a,b,c\rangle \star b.$$
Hence for any $a,b,c\in A,$
$$\langle a,b\star b,c\rangle=2 b\star \langle a,b,c\rangle.$$

Prove now  that, if $A$ is assosymmetric, then 
$$glennie(a,b,c)=shest(a,b,c\star c)-2\,c\star shest(a,b,c)=0$$
is polynomial identity for plus-assosymmetric algebras. 
For $p\ne 2,3,$ we can apply results of Section \ref{Associators}. 
 
Let $p\ne 2,3.$   By Lemma \ref{five}, $D(a,b,c)=shest(a,b,c)$ is well defined on plus-assosymmetric algebras and the map
$c\mapsto D(a,b,c) $ is a derivation of assosymmeric algebras. In particular, it is  a derivation of plus-assosymmetric algebras. 
The condition $shest(a,b,c\star c)=0$ is equivalent to Leibniz condition  for this derivation.  
As was shown by Shestakov,  the identity $shest=0$ is equivalent to  Glennie identity \cite{Glennie}. 
Suppose that  the identity $shest=0$ is a consequence of commutativity and Lie-triple identities. Since Lie-triple identity is a consequence of 
Jordan identity, this means that $shest=0$ will be identity for Jordan algebars also. It contradicts to Glennie's result that Glennie identity is exceptional. 
Therefore, commutative, Jordan and Glennie identities are independent non only in a class of plus-associative algebras, but also in a class of plus-assosymmetric algebras. 

For $p=3$ one can check by computer program, say by albert \cite{albert}, that Glennie identity holds for plus-assosymmetric algebras also.

\section{Identities of degree 4 for commutative algebras}

In this section we consider non-associative but {\it commutative} polynomials with variables $t_1,t_2,\cdots .$ 
Let $wjor$ be multilinear Jordan polynomial
$$wjor(t_1,t_2,t_3,t_4)=(t_2,t_1,t_3t_4)+(t_3,t_1,t_4t_2)+(t_4,t_1,t_2t_3)$$
and 
$$jor_2(t_1,t_2,t_3,t_4)=$$
$$
wjor(t_1,t_2,t_3,t_4)-wjor(t_2,t_1,t_3,t_4)+wjor(t_3,t_1,t_2,t_4)-wjor(t_4,t_1,t_2,t_3).$$
Define multilinear commutative non-associative polynomial $jor_1$ by 
$$jor_1(t_1,t_2,t_3,t_4)=[l_{t_1},l_{t_2}](t_3t_4)-t_3([l_{t_1},l_{t_2}](t_4))-([l_{t_1},l_{t_2}](t_3))t_4=$$
$$t_1(t_2(t_3t_4))-t_3(t_1(t_2t_4))-(t_1(t_2t_3))t_4
-t_2(t_1(t_3t_4))-t_3(t_2(t_1t_4))-(t_2(t_1t_3))t_4.$$
Note that $jor_1(t_1,t_2,t_3,t_4)$ is skew-symmetric by variables $t_1,t_2$ and symmetric by variables $t_3,t_4.$ 

\begin{lm} 
\label{1111}
If $A$ is commutative algebra, then for any $a,b,c\in A,$  
$$(a,b,c)=[l_a,l_c](b).$$
In particular, for any $a,b,c\in A,$ 
$$(a,b,c)+(c,b,a)=0.$$
\end{lm}

\sudda{{\bf Proof.} 
$$(a,b,c)=a(bc)-(ab)c=a(cb)-c(ab)=l_al_c(b)-l_cl_a(b)=[l_a,l_c](b).$$}

\begin{lm}\label{1122}
$$jor_1(t_1,t_2,t_3,t_4)=(t_1,t_3t_4,t_2)-t_3(t_1,t_4,t_2)-t_4(t_1,t_3,t_2).$$
\end{lm}

{\bf Proof.} Easy calculations based on  Lemma \ref{1111}. 

\begin{lm}\label{1112}
$jor_1(t_1,t_2,t_3,t_4)=-wjor(t_1,t_2,t_3,t_4)+wjor(t_2,t_1,t_3,t_4).$
\end{lm}

{\bf Proof.}  We have 
$$wjor(t_1,t_2,t_3,t_4)-wjor(t_2,t_1,t_3,t_4)=$$
$$(t_2,t_1,t_3t_4)+(t_3,t_1,t_4t_2)+(t_4,t_1,t_2t_3)-$$
$$(t_1,t_2,t_3t_4)-(t_3,t_2,t_4t_1)-(t_4,t_2,t_1t_3)=$$
$$t_2(t_1(t_3t_4))-(t_2t_1)(t_3t_4)+t_3(t_1(t_4t_2))-(t_3t_1)(t_4t_2)
+t_4(t_1(t_2t_3))-(t_4t_1)(t_2t_3)$$
$$-t_1(t_2(t_3t_4))+(t_1t_2)(t_3t_4)-t_3(t_2(t_4t_1))+(t_3t_2)(t_4t_1)
-t_4(t_2(t_1t_3))+(t_4t_2)(t_1t_3)=$$
$$t_2(t_1(t_3t_4))+t_3(t_1(t_4t_2))+t_4(t_1(t_2t_3))$$
$$-t_1(t_2(t_3t_4))-t_3(t_2(t_4t_1))-t_4(t_2(t_1t_3))=$$
$$-t_1(t_2(t_3t_4))+t_3(t_1(t_2t_4))+(t_1(t_2t_3))t_4$$
$$t_2(t_1(t_3t_4))-t_3(t_2(t_1t_4))-(t_2(t_1t_3))t_4=$$
$$-jor_1(t_1,t_2,t_3,t_4).$$

\begin{lm} \label{2222}
$$jor_2(t_1,t_2,t_3,t_4)-jor_2(t_2,t_1,t_3,t_4)=-2\,jor_1(t_1,t_2,t_3,t_4),$$
$$jor_2(t_1,t_2,t_3,t_4)-jor_2(t_1,t_2,t_4,t_3)= -2\,jor_1(t_3,t_4,t_1,t_2) ,$$
$$jor_2(t_1,t_2,t_3,t_4)+jor_2(t_2,t_1,t_3,t_4)=-2\,jor_1(t_3,t_4,t_1,t_2),$$
$$jor_2(t_1,t_2,t_3,t_4)+jor_2(t_1,t_2,t_4,t_3)=-2\, jor_1(t_1,t_2,t_3,t_4).$$
In particular,
\begin{equation}\label{xxx}
jor_2(t_1,t_2,t_3,t_4)+jor_2(t_2,t_1,t_4,t_3)=0.
\end{equation}
\end{lm}

{\bf Proof.} 
Since variables $t_i$ are commuting,  by Lemma \ref{1111} 
$$jor_2(t_1,t_2,t_3,t_4)-jor_2(t_2,t_1,t_3,t_4)=
2(wjor(t_1,t_2,t_3,t_4)-wjor(t_2,t_1,t_3,t_4))$$
Therefore by Lemma \ref{1112},
$$jor_2(t_1,t_2,t_3,t_4)-jor_2(t_2,t_1,t_3,t_4)
=-2\,jor_1(t_1,t_2,t_3,t_4).$$
Further, by Lemma \ref{1112},
$$jor_2(t_1,t_2,t_3,t_4)-jor_2(t_1,t_2,t_4,t_3)=$$
$$2\{wjor(t_3,t_1,t_2,t_4)-wjor(t_4,t_1,t_2,t_3)\}=$$
$$2\{wjor(t_3,t_4,t_1,t_2)-wjor(t_4,t_3,t_1,t_2)\}=$$
$$ -2\,jor_1(t_3,t_4,t_1,t_2).$$
Similarly, by Lemma \ref{1111} and Lemma \ref{1112},
$$jor_2(t_1,t_2,t_3,t_4)+jor_2(t_2,t_1,t_3,t_4)=$$
$$2\{wjor(t_3,t_1,t_2,t_4)-wjor(t_4,t_1,t_2,t_3)\}=$$
$$2\{wjor(t_3,t_4,t_1,t_2)-wjor(t_4,t_3,t_1,t_2)\}=$$
$$-2\,jor_1(t_3,t_4,t_1,t_2),$$
and 
$$jor_2(t_1,t_2,t_3,t_4)+jor_2(t_1,t_2,t_4,t_3)=$$
$$2\{wjor(t_1,t_2,t_3,t_4)-wjor(t_2,t_1,t_3,t_4)\}=$$
$$-2\,jor_1(t_1,t_2,t_3,t_4).$$

Since
$$jor_2(t_1,t_2,t_3,t_4)-jor_2(t_1,t_2,t_4,t_3)= -2\,jor_1(t_3,t_4,t_1,t_2) =$$ $$
jor_2(t_1,t_2,t_3,t_4)+jor_2(t_2,t_1,t_3,t_4),$$
we have  
$$jor_2(t_1,t_2,t_4,t_3)+jor_2(t_2,t_1,t_3,t_4)=0.$$
$\square$

\begin{lm}\label{3333}
$$jor_2(t_1,t_2,t_3,t_4)=- jor_1(t_1,t_2,t_3,t_4)-jor_1(t_3,t_4,t_1,t_2).$$
\end{lm}

{\bf Proof.}By Lemma \ref{2222}
$$jor_2(t_1,t_2,t_3,t_4)-jor_2(t_2,t_1,t_3,t_4)=-2\,jor_1(t_1,t_2,t_3,t_4),$$
$$jor_2(t_1,t_2,t_3,t_4)-jor_2(t_1,t_2,t_4,t_3)= -2\,jor_1(t_3,t_4,t_1,t_2) ,$$
 Add these two relations. By Lemma \ref{2222}, relation (\ref{xxx}), we receive  that 
$$jor_2(t_1,t_2,t_3,t_4)=- jor_1(t_1,t_2,t_3,t_4)-jor_1(t_3,t_4,t_1,t_2).$$
$\square$

\begin{lm} \label{Arman} Let $p\ne 2$ and  $A$ be a commutative algebra.
Then the following conditions are equivalent
\begin{equation}
\label{11}
[l_a,l_b]\in Der A
\end{equation}
\begin{equation}
\label{22}
\sudda{jor_1(a,b,c,d)=0}
a(b(cd))-(a(bc))d-c(a(bd))=
b(a(cd))-(b(ac))d-c(b(ad))
\end{equation}
\begin{equation}
\label{33}
\sudda{jor_1(a,b,c,c)=0}
 a(b(cc))-2(a(bc))c=b(a(cc))-2c(b(ac))
\end{equation}
\begin{equation}
\label{44}
(a,cc,b)=2c(a,c,b)
\end{equation}
\begin{equation}
\label{55}
(a,bc,d)=b(a,c,d)+c(a,b,d)
\end{equation}
\begin{equation}
\label{66}
jor_2(a,b,c,d)=0
\end{equation}
\begin{equation}
\label{99}
[[l_a,l_c],l_b]=l_{(a,b,c)}
\end{equation}
\begin{equation}
\label{100}
wjor(a,b,c,d)=wjor(b,a,c,d)
\end{equation}
Any of these identities implies the identity
\begin{equation}
\label{77}
2((ba)a)a+b((aa)a)=3(b(aa))a 
\end{equation}
\end{lm}

{\bf Proof.} The Leibniz rule for a derivation $[l_a,l_b]$ can be written as 
$$[l_a,l_b](cd)=([l_a,l_b]c)d+c([l_a,l_b]d)$$
Therefore  (\ref{11}) is equivalent to (\ref{22}). 
Substitution in (\ref{22}) $c=d$ gives us (\ref{33}). Conversely,  polarization of (\ref{33}) 
gives us (\ref{22}). 
Rewrite (\ref{33}) by the following way,
$$a(b(cc))-b(a(cc))=2c\{ a(bc)-b(ac)\}.$$
For commutative algebra this condition is equivalent to (\ref{44}). 
The identity (\ref{55}) is a polarization of (\ref{44}),
$$(a,(b+c)^2,d)-2(b+c)(a,b+c,d)-(a,b^2,d)+2b(a,b,d)-(a,c^2,d)+2c(a,c,d)=$$
$$(a,bc,d)+(a,cb,d)-2b(a,c,d)-2c(a,b,d)=$$
$$2\{(a,bc,d)-b(a,c,d)-c(a,b,d)\}.$$
So, (\ref{44}) implies (\ref{55}). If $p\ne 2,$  let us substitute $b=c$ in (\ref{55}). We obtain  
(\ref{44}).  

By Lemma \ref{1112}, identities $jor_1=0$ and $jor_2=0$ are equivalent. 
Therefore, by Lemma \ref{1122} identities (\ref{55}) and (\ref{66}) are equivalent. 
Further
$$jor_2(a,b,c,c)=(a,cc,b)-2c(a,c,b).$$
So, (\ref{44}), (\ref{55}), (\ref{66}) are equivalent. 
We have 
$$([[l_a,l_c],l_b]-l_{(a,b,c)})d=$$
$$a(c(bd))-c(a(bd))-b(a(cd))+b(c(ad))-(a,b,c)d=$$
$$(a,bd,c)-b(a,d,c)-(a,b,c)d.$$
Therefore (\ref{99}) is equivalent to (\ref{55}). 
By Lemma  \ref{1112} conditions (\ref{22}) and  (\ref{100}) are equivalent. 
Take in (\ref{66}) $c=d=a,$
$$jor_2(a,b,a,a)=$$
$$(a,aa,b)-2 a(a,a,b)+(a,ab,a)-a(a,b,a)-b(a,a,a)=$$
$$a((aa)b)-(a(aa))b-2 a(a(ab))+2a((aa)b)=$$
$$3a(b(aa))-b(a(aa))-2a(a(ab)).$$
So, (\ref{66}) implies (\ref{77}).

{\bf Remark.} If $p>3$ then (\ref{77}) and any of identities  {\rm (\ref{11}), (\ref{22}), (\ref{33}), (\ref{44}), (\ref{55}), (\ref{66}) } and {\rm (\ref{99})} are equivalent.  

{\bf Remark.}  In Lemma \ref{Arman} we assume that the commutativity identity $t_1t_2=t_2t_1$ is given. If we omit the commutativity identity, then identities (\ref{11}) --  (\ref{100}) are not equivalent. 
Namely,  (\ref{11}) and (\ref{44}) are equivalent,  (\ref{22}), (\ref{33}), (\ref{55}), (\ref{66}) are equivalent,  (\ref{11}) and (\ref{44}) imply (\ref{22}), (\ref{33}), (\ref{55}), (\ref{66}) and 
(\ref{22}), (\ref{33}), (\ref{55}), (\ref{66}) do not imply (\ref{11}) and (\ref{44}).

\section{Proof of Theorem \ref{q^2=1}.}

Consider identities for the class ${\mathcal A}ss^{+}.$ It has identity of degree 2: commutativity identity. No identity of degree $3$ that does not follow from commutativity rule. In degree $4$ 
the following commutative non-associative polynomials gives us identities 
 that are not consequences of commutativity identity 
(see \cite{Jacobson}, chapter 1.1, p.5-6):
$$g_{[4]}^{(1)}(t_1)=(t_1,t_1,t_1^2),$$

$$ g_{[3,1]}^{(1)}(t_1,t_2)=(t_1,t_2,t_1^2),$$
$$g_{[3,1]}^{(2)}(t_1,t_2)=t_2(t_1t_1^2)+2t_1(t_1(t_1t_2))-3t_1(t_2t_1^2),$$

$$g_{[2,2]}^{(1)}(t_1,t_2)=t_1^2t_2^2-t_1(t_1t_2^2)-2t_2(t_1(t_1t_2))+2(t_1t_2)(t_1t_2),$$

$$g_{[2,1,1]}^{(1)}(t_1,t_2,t_3)=(t_1,t_1,t_2t_3)+(t_2,t_1,t_3t_1)+(t_3,t_1,t_1t_2),$$
$$g_{[2,1,1]}^{(2)}(t_1,t_2,t_3)=2(t_1,t_2,t_1t_3)+(t_3,t_2,t_1^2),$$

$$g_{[1,1,1,1]}^{(1)}(t_1,t_2,t_3,t_4)=(t_2,t_1,t_2t_3)+(t_3,t_1,t_1t_4)+(t_4,t_1,t_2t_3).$$
All polynomials  of a form $g_{\alpha}^{(i)}$    are homogeneous of type $\alpha.$ This means that lower index $\alpha$ corresponds to the  type of identity, 
i.e., if $\alpha=[\alpha_1,\ldots,  \alpha_s,\ldots ],$ then $t_s$ in each monomial of $g_{\alpha}^{(i)}$  enters $\alpha_s$ times. Such $\alpha_s$ is called multiplicity of $t_s.$
Note that permutation of indices $t_s$ with equal multiplicities in $g_{\alpha}^{(i)}$ induces a consequence of the identity $g_{\alpha}^{(i)}$ of type $\alpha.$ Therefore,  any consequence of these identities in degree $4$  can be presented as a linear combination of polynomials with given type where variables with equal multiplicity are permuted. 

Any associative algebra is assosymmetric. Therefore, any identity of type $\alpha$  of degree 4 for plus-assosymmetric  algebras is a consequence of identities $g_\alpha^{(i)}.$ So,  polynomials of the following form should be tested for an identity of plus-assosymmetric  algebras
$$\alpha=[4],  \qquad f_{[4]}=g_{[4]}^{(1)},$$
$$\alpha =[3,1], \qquad f_{[3,1]}^{\mu_1,\mu_2}=\mu_1 g_{[3,1]}^{(1)}+\mu_2 g_{[3,1]}^{(2)},\quad \mu_i\in K,$$
$$\alpha=[2,2], \qquad f_{[2,2]}^{\mu_1,\mu_2}(t_1,t_2)=\mu_1g_{[2,2]}^{(1)}(t_1,t_2)+\mu_2g_{[2,2]}^{(1)}(t_2,t_1),\quad \mu_i\in K.$$
Since $g_{[2,2]}^{(1)}(t_1,t_2,t_3)=g_{[2,2]}^{(1)}(t_1,t_3,t_1),$ 
in case of  $\alpha=[2,1,1],$ as a general form of a commutative polynomial tested for identity of plus-assosymmetric  algebras we can get 
$$ f_{[2,1,1]}^{\mu_1,\mu_2,\mu_3}(t_1,t_2,t_3)=$$
$$\mu_1 g_{[2,1,1]}^{(1)}(t_1,t_2,t_3)+\mu_2 g_{[2,1,1]}^{(2)}(t_1,t_2,t_3)+\mu_3 g_{[2,1,1]}^{(2)}(t_1,t_3,t_2), \quad \mu_i\in K.$$
Recall that  $wjor(t_1,t_2,t_3,t_4)$ are symmetric by permutations of indices $t_2,t_3,t_4.$ 
Therefore as a general form of a commutative polynomial tested for identity of plus-assosymmetric  algebras in the case $\alpha=[1,1,1,1]$  we can take 
 $$ f_{[1,1,1,1]}^{\mu_1,\mu_2,\mu_3,\mu_4}(t_1,t_2,t_3,t_4)=$$
 $$\mu_1 g_{[1,1,1,1]}^{(1)}(t_1,t_2,t_3,t_4)+
 \mu_2 g_{[1,1,1,1]}^{(1)}(t_2,t_1,t_3,t_4)+\mu_3 g_{[1,1,1,1]}^{(1)}(t_3,t_1,t_2,t_4)+$$
 $$\mu_4 g_{[1,1,1,1]}^{(1)}(t_4,t_1,t_2,t_3), \quad \mu_i\in K.$$ 

By theorem 1 of \cite{Hentzel} free assosymmetric  algebra in type $\alpha$ of  degree $4$  have the following base and dimensions

$$\begin{array}{|c|l|c|}
\hline
\alpha&base&dim\\
\hline
[4]&\{((aa)a)a, (aa)(aa), (a(aa))a\}&3\\
\hline
[3,1]&\{(aa)(ab),  (b(aa))a,    ((ba)a)a, ((ab)a)a,     ((aa)b)a,     (a(aa))b, &7\\  
& ((aa)a)b\}&\\
\hline
[2,2]&\{
\{(aa)(bb),
  (b(ab))a,
((bb)a)a,
((ba)b)a,
((ab)b)a,
(b(aa))b,&9\\
&((ba)a)b,
 ((ab)a)b,
((aa)b)b\}&\\
\hline
[2,1,1]&
\{  (aa)(bc),
(c(ab))a,
((cb)a)a,
 ((bc)a)a,
 ((ca)b)a,&16\\
&((ac)b)a,
((ba)c)a,
  ((ab)c)a,
 (c(aa))b,
((ca)a)b,
((ac)a)b,&\\
&((aa)c)b,
(b(aa))c,
 ((ba)a)c,
((ab)a)c,
  ((aa)b)c\}&\\
\hline
[1,1,1,1]&\mbox{too big to be presented here}&29\\
\hline
\end{array}$$

\medskip

Let us substitute in  polynomials $f_{\alpha}$ instead of parameters $t_i$ elements of free assosymmetric algebras and calculate its value in terms of assosymmetric multiplication. 

We have 
$$f_{[4]}(a)=<a,a,\{a,a\}>=$$
$$\{a,\{a,\{a,a\}\}\}-\{\{a,a\},\{a,a\}\}=$$
$$2(a(a(aa))+2a((aa)a)+2(a(aa))a+2((aa)a)a-8(aa)(aa)=$$
$$2(a,a,aa)-2(aa,a,a)+2a((aa)a)+2(a(aa))a-4(aa)(aa)=$$
$$2(a,aa,a)+4(a(aa))a-4(aa)(aa)=$$
$$2(aa,a,a)+4(a(aa))a-4(aa)(aa)=$$
$$2(aa)(aa)-2((aa)a)a+4(a(aa))a-4(aa)(aa)=$$
$$-2((aa)a)a+4(a(aa))a-2(aa)(aa).$$
Since elements $((aa)a)a, (a(aa))a$ and $(aa)(aa)$ are base elements, this means that 
$f_{[4]}(a)\ne 0.$
So, plus-assosymmetric  algebras have no identity of type [4].

Similar calculations show that $f_{[3,1]}=0$ is identity if $\mu_1=0.$ So, in type [3,1] plus-assosymmetric  algebras have an identity $g_{[3,1]}^{(2)}=0.$

Consider  type [2,2] case. 
Calculations show that 
$$f_{[2,2]}^{\mu_1,\mu_2}(t_1,t_2)=$$
$$6(\mu_1+\mu_2) \{(a a)(b, b) - 2(b(ab)) a -   ((a a) b) b + 2 ((ba) b) a\}.$$
So, $f_{[2,2]}^{\mu_1,\mu_2}$ is identity for plus-assosymmetric  algebras if $\mu_1+\mu_2=0,$ $p\ne 2,3,$ 
and,  $h_{[2,2]}(t_1,t_2)=0$ is identity for plus-assosymmetric  algebras, where 
$h_{[2,2]}=f_{[2,2]}^{-1,1}.$ Note that 
$$h_{[2,2]}(t_1,t_2)=t_1(t_1t_2^2)-t_2(t_2t_1^2)-2t_1(t_2(t_1t_2))+2t_2(t_1(t_1t_2))$$
is identity for plus-assosymmetric  algebras. 

In the case  of type [2,1,1], we have 
$$f_{[2,1,1]}^{\mu_1,\mu_2,\mu_3}(a,b,c)=
-6 (\mu_1 + \mu_2 + \mu_3) \{(aa)(bc) - 2(c(ab)) a - 
   ((aa)b)c+ 2 ((c a) b) a\}.$$
Therefore,  $h_{[2,1,1]}^{(1)}=0$ and $h_{[2,1,1]}^{(2)}=0$ 
are identities for plus-assosymmetric  algebras, where 
$$h_{[2,1,1]}^{(1)}=f_{[2,1,1]}^{1,-1,0},\qquad h_{[2,1,1]}^{(2)}=f_{[2,1,1]}^{0,1,-1},$$
if $p\ne 2,3.$ 
Note that 
$$h_{[2,1,1]}^{(1)}(t_1,t_2,t_3)=(t_1,t_1,t_2t_3)+(t_2,t_1,t_3t_1)+(t_3,t_1,t_1t_2)
-2(t_1,t_2,t_1t_3)-(t_3,t_2,t_1^2),$$
$$h_{[2,1,1]}^{(2)}(t_1,t_2,t_3)=$$ 
$$2(t_1,t_2,t_1t_3)-2(t_1,t_3,t_1t_2)+(t_3,t_2,t_1^2)-(t_2,t_3,t_1^2).$$

Now consider the case [1,1,1,1].  By Lemma \ref{4} the polynomial $wjor$ can not give an identity on ${\mathcal A}ssym^+,$ if $p\ne 2,3.$ 
For any assosymmetric algebra $A$ and for any its four elements $a,b,c,d\in A$ by Lemma~\ref{5}, 
$$f_{[1,1,1,1]}^{\mu_1,\mu_2,\mu_3,\mu_4}(a,b,c,d)=(\mu_1+\mu_2+\mu_3+\mu_4)wjor(a,b,c,d).$$
Thus, $f_{[1,1,1,1]}^{\mu_1,\mu_2,\mu_3,\mu_4}=0$ is identity on ${\mathcal A}ssym^+$ if and only if  
$$\mu_1+\mu_2+\mu_3+\mu_4=0$$
Therefore
$$\mu_4=-\mu_1-\mu_2-\mu_3,$$
and 
$$f_{[1,1,1,1]}^{\mu_1,\mu_2,\mu_3,\mu_4}(t_1,t_2,t_3,t_4)=$$
$$\mu_1 (wjor(t_1,t_2,t_3,t_4)- wjor(t_4,t_1,t_2,t_3))+$$
$$\mu_2( wjor(t_2,t_1,t_3,t_4)- wjor(t_4,t_1,t_2,t_3))+$$
$$\mu_3 (wjor(t_3,t_1,t_2,t_4) - wjor(t_4,t_1,t_2,t_3)).$$
In other  words, by Lemma \ref{1112}, 
$$-f_{[1,1,1,1]}^{\mu_1,\mu_2,\mu_3,\mu_4}(t_1,t_2,t_3,t_4)=$$
$$\mu_1\, jor_1(t_1,t_4,t_2,t_3)+\mu_2\,jor_1(t_2,t_4,t_1,t_3)+
\mu_3\,jor_1(t_3,t_4,t_1,t_2).$$
By Lemma \ref{5}  $jor_1=0$ is identity on ${\mathcal A}ssym^+.$ 
So,  $jor_1=0$ is   identity for ${\mathcal A}ssym^+.$

It remains to  prove that all identities appeared for types [2,2] and [2,1,1] are consequences of the identity $jor_1=0$ if $p\ne 2,3.$ 

It is easy too see that 
$$h_{[2,2]}(t_1,t_2)=$$
$$(t_1,t_1,t_2^2)-(t_2,t_2,t_1^2)-2(t_1,t_2,t_1t_2)+2(t_2,t_1,t_1t_2)=$$
$$t_1(t_1t_2^2)-t_2(t_2t_1^2)-2\{t_1(t_2(t_1t_2))-t_2(t_1(t_1t_2))\}= $$
$$-t_1(t_2,t_2,t_1)+t_2(t_1,t_1,t_2)-t_1(t_2(t_1t_2))+t_2(t_1(t_1t_2)).$$
Therefore  by Lemma \ref{1111},
$$h_{[2,2]}(t_1,t_2)=$$
$$t_1(t_1,t_2,t_2)+t_2(t_1,t_1,t_2)-t_1((t_1t_2)t_2)+(t_1(t_1t_2))t_2=$$
$$t_1(t_1,t_2,t_2)+t_2(t_1,t_1,t_2)-(t_1,t_1t_2,t_2).$$
Hence the identity $h_{[2,2]}=0$ is a consequence of the identity 
$$(t_1,t_2t_3,t_4)-t_1(t_1,t_3,t_4)-(t_1,t_2,t_4)t_3=0.$$
By Lemma \ref{Arman}, (\ref{55}),  this means that the identity $h_{[2,2]}=0$ is a consequence of 
identity $jor_1=0.$
By Lemma \ref{1112}
$$h_{[2,1,1]}^{(1)}=  wjor(t_1,t_2,t_3,t_1)-wjor(t_2,t_1,t_3,t_1)=-
jor_1(t_1,t_2,t_3,t_1),$$ 
$$h_{[2,1,1]}^{(2)}=wjor(t_2,t_3,t_1,t_1)-wjor(t_3,t_2,t_1,t_1)=-jor_1(t_2,t_3,t_1,t_1).$$
So, $jor_1=0$ is 
a minimal identity for ${\mathcal A}ssym^+$ that does not follow from commutativity identity if $p=char\,K\ne 2,3.$ 

Now consider the case $p=3.$ By Lemma \ref{4} $wjor=0$ is an identity for plus-assosymmetric  algebras. As we have checked above $g_{[2,2]}^{(1)}=0$ and $g_{[2,1,1]}^{(1)}=g_{[2,1,1]}^{(2)}=0$ are identities for $p=3.$ 
We have 
$$g_{[2,2]}^{(1)}(t_1,t_2)=-(t_1,t_1,t_2^2)-2(t_2,t_1,t_1t_2)=-wjor(t_1,t_1,t_2,t_2),$$
$$g_{[2,1,1]}^{(1)}(t_1,t_2,t_3)=wjor(t_1,t_1,t_2,t_3),$$
$$g_{[2,1,1]}^{(2)}(t_1,t_2,t_3)=2(t_1,t_2,t_1t_3)+(t_3,t_2,t_1^2)=wjor(t_2,t_3,t_1,t_1).$$
So, $wjor=0$ is a minimal identity for ${\mathcal A}ssym^+,$  $p=3.$

Since by Lemma \ref{Arman} Lie triple identity and the identity $jor_1=0$ are equivalent, Theorem \ref{q^2=1} is proved completely. 

\section{Additional remarks}

{\bf Remark 1} concerns Theorem \ref{q^2=1} in case $p=2.$  If $p=2,$ then $\mathcal{A}ssym^{(+)}=\mathcal{A}ssym^{(-)}.$ Therefore,  plus-assosymmetric algebras are Lie, and 
all identities for ${\mathcal A}ssym^{(+)}$ follows from commutativity and Jacobi identities. In particular, Jordan identity and Glennie identity 
are consequences of commutativity one and Jacobi identity.

{\bf Remark 2.} Recall that an Jordan algebra is called {\it special} Jordan, if it is isomorphic to a subalgebra of algebra $A^{(+)}$ for some associative algebra $A.$  Well known that Jordan algebra of $3\times 3$ hermitian matrices over octonians $M_3^8$ is not special, and Glennie identity is an example of special Jordan identity.  Say that a Lie triple algebra is {\it special} if it is isomorphic to a subalgbera of algebra $A^{(+)}$ for some Lie triple algebra $A.$  By Theorem \ref{main1} $M_3^8$ as Lie triple algebra is exeptional and Glennie identity is special Lie triple identity.

{\bf Remark 3.} Non-commutative Jordan algebras are defined as flexible algebras with identity $jor=0.$ The flexible identity appears as a consequence of Jordan identity if one assumes existence of unit. For Lie triple case, supplementing of unite gives nothing, and one can define non-commutative versions of Lie triple algebras in two ways. Let 
$$assder(t_1,t_2,t_3,t_4)=(t_1,t_2t_3,t_4)-t_2(t_1,t_3,t_4)-(t_1,t_2,t_4)t_3.$$
An algebra with identity $assder=0$ is called {\it associator is derivation} algebra. The identity $lietriple=0$ is a consequence of the identity $assder=0.$
One can consider non-commutative version of Lie triple algebras as algebras with identity $assder=0$ or one can understand non-commutative Lie triple algebras as algebras with identity $lietriple=0.$  We think that second version is perefereable. It contains a class of non-commutative Jordan algebras. 

For an algebra $A=(A,\cdot),$  let us define its $q$-commutator by 
$a\cdot_q b=a\cdot b+q\,b\cdot a.$ Quasi-associative algebras are defined as associative algebras under $q$-commutator.
Similarly, let us define quasi-assosymmetric algebras as assosymmetric algebras under $q$-commutator. 

 Let $\sigma_q:K\langle t_1,\ldots,t_k\rangle \rightarrow K\langle t_1,\ldots,t_k\rangle$ is endomorphism of free magmatic algebra that corresponds to a magmatic polynomial  a  polynomial where magmatic products are changed my their $q$-commutators. For example, if 
$f=(t_2\cdot t_3)\cdot t_1$ then $\sigma_q f=(t_2\cdot_q t_3)\cdot_q t_1=(t_2\cdot t_3)\cdot t_1+q\,(t_3\cdot t_2)\cdot t_1+q\,t_1\cdot (t_2\cdot t_3)+q^2t_1\cdot(t_3\cdot t_2).$
 In \cite{qleibniz} Theorems 2.2 and 2.3 state that if $q^2\ne 1,$ and  $\mathcal C=Var(f_1,\ldots,f_s)$ is a variety, then
 \begin{itemize}
 \item   ${\mathcal C}^{(q)}$ is a variety also, and 
 \item ${\mathcal C}^{(q)}=Var(\sigma_{-q}f_1, \ldots, \sigma_{-q}f_s).$ 
 \end{itemize}
By this result, quasi-assosymmetric algebras forms variety if $q^2\ne 1,$  and the variety ${\mathcal A}ssym^{(q)}$ is generated by identites $lsym^{(q)}=0$ and 
$r\!sym^{(q)}=0,$ where 
$$lsym^{(q)}(t_1,t_2,t_3)=\sigma_{-q}lsym(t_1,t_2,t_3)=$$ 
$$t_1(t_2 t_3) - q t_1(t_3 t_2)
-t_2(t_1 t_3) + q t_2(t_3 t_1)
+$$ $$q (q + 1) t_3(t_1 t_2)
-q (q + 1) t_3(t_2 t_1)
-(1 + q) (t_1 t_2) t_3
+$$ $$(q + 1) (t_2 t_1) t_3
+q (t_1 t_3) t_2- q^2 (t_3 t_1) t_2
-q (t_2 t_3) t_1 + q^2 (t_3 t_2) t_1$$
and 
$$r\!sym^{(q)}=\sigma_{-q}rsym(t_1,t_2,t_3)=$$ 
$$
(q + 1) t_1(t_2 t_3)
-(q + 1) t_1(t_3 t_2)
-q t_2(t_1 t_3) + q^2 t_2(t_3 t_1)
+q t_3(t_1 t_2) - q^2 t_3(t_2 t_1)
$$
$$
-(t_1 t_2) t_3 + (t_1 t_3) t_2
+q (t_2 t_1) t_3
-q (t_3 t_1) t_2
-q (q + 1) (t_2 t_3) t_1
+q (q + 1) (t_3 t_2) t_1.$$
Note that the identity $assder=0$ is a consequence of identities $lsym^{(q)}=0, r\!sym^{(q)}=0.$ Therefore, quasi-assosymmetric algebras gives us another examples of non-commutative Lie triple algebras. 

{\bf Remark 4.}  Assosymmetric operad is not Koszul.  Let us give sketch of proof. 
Dual operad to assosymmetric operad is generated by identities 
$$[a,b]c+[b,c]a+[c,a]b=0,$$
$$(a,b,c)=0.$$
Let $d_n^!$ are dimensions of multilinear part of free algebra with such identities. Then 
$$d_1^!=1, d_2^!=2, d_3^!=5, d_4^!=9, d_5^!=9, d_6^!=11, d_7^!=13.$$
Generating functions of assosymetric and dual assostymetric operads looks like 
$$G_{assym}(x)=-x+2x^2/2!-7x^3/3!+29 x^4/4!-136x^5/5!+O(x)^6,$$
$$G_{assym}^!(x)=-x+2x^2-5x^3/3!+9x^4/4!-9x^5/5!+O(x)^6.$$
and 
$$G_{assym}(G_{assym}^!(x))=x+3x^5/8+O(x)^6\ne x.$$
So, by Koszulity criterium ( \cite{GiK} Proposition 4.14(b) )
assosymmetric operad is not Koszul.

\end{document}